\documentclass[letterpaper, 10 pt, conference]{ieeeconf}%[twocolumn]{IEEEtran}
\usepackage[T1]{fontenc}% optional T1 font encoding

\IEEEoverridecommandlockouts                          
\usepackage{amsmath} % assumes amsmath package installed
\usepackage{amssymb}  % assumes amsmath package installed
\usepackage{xcolor}
\usepackage{caption}
\usepackage{graphicx}
\usepackage{mathtools}
\usepackage{caption}
\usepackage{floatrow}
\usepackage{tikz}

\usepackage{algorithmicx}
    \usepackage[ruled]{algorithm}
    \usepackage{algpseudocode}

\usepackage{pgfplots}
\usepackage{xcolor}
\usetikzlibrary{matrix,arrows,calc,positioning,shapes,decorations.pathreplacing}
\usepackage{graphicx}

\usepackage{amsmath} % assumes amsmath package installed

\usepackage{enumerate}
\usepackage[all,tips]{xy}
\SelectTips{cm}{11}

\newtheorem{definition}{Definition}
\newtheorem{remark}{Remark}

\newcommand{\g}{\mathfrak{g}}

\newcommand{\Vezero}{U_{i,0}^{\text{ext}}}
\newcommand{\R}{\mathbb{R}}
\newcommand{\Ad}{\text{Ad}}
\newcommand{\ad}{\text{ad}}
\newcommand{\se}{\mathfrak{se}(2)}

\DeclareMathOperator{\tr}{tr}

\begin{document}

\title{
Optimal Control with Broken Symmetry of\\Multi-Agent Systems on Lie Groups}

\author{Efstratios Stratoglou, Leonardo Colombo and Tomoki Ohsawa 
\thanks{E. Stratoglou (ef.stratoglou@alumnos.upm.es) is with Universidad Polit\'ecnica de Madrid (UPM), 28006 Madrid, Spain. L. Colombo (leonardo.colombo@car.upm-csic.es) is with Centre for Automation and Robotics (CSIC-UPM), Ctra. M300 Campo Real, Km 0,200, Arganda
del Rey - 28500 Madrid, Spain. T. Ohsawa (tomoki@utdallas.edu) is with Department of Mathematical Sciences, The University of Texas at Dallas, 800 W Campbell Rd,
Richardson, TX 75080-3021. This work was partially supported by NSF grant CMMI-1824798.}}%

\maketitle

\begin{abstract}
    
In this paper we study reduction by symmetry for optimality conditions in  optimal control problems of left-invariant affine multi-agent control systems, with partial symmetry breaking cost functions. Our approach emphasizes the role of variational principles. Specifically, we recast the optimal control problem as a constrained variational problem with a partial symmetry breaking Hamiltonian and obtain the reduced optimality conditions from a reduced variational principle via Pontryagin Maximum Principle. We apply the results to a collision avoidance problem for multiple unicycles in the presence of an obstacle.\end{abstract}

%\begin{IEEEkeywords} Symmetry reduction, multi-agent systems, optimal control, collision and obstacle avoidance. \end{IEEEkeywords}

\IEEEpeerreviewmaketitle

\section{Introduction}

 Lie groups symmetries appear naturally in many control systems problems \cite{B}, 
 \cite{borum},  \cite{grizzle}. Examples of invariant control problems on Lie groups include motion planning for underwater vehicles \cite{Leonard1}, conflict resolution in differential games \cite{tomlin}, collective motion in biological models \cite{JK2}, controlled Lagrangians \cite{BlControlled,contreras}, and coordination of multi-agent systems \cite{bonnabel}.
 
 The goal of this work is to study symmetry reduction in optimal control problems for left-invariant affine multi-agent control systems on Lie groups with partially broken symmetries. More specifically, cost functions that break some of the symmetries but not all. Symmetry breaking is a common phenomena in several physical contexts \cite{GBT}. The simplest example is the spinning top dynamics (the motion of a rigid body with a fixed point under a gravitational field), where due to the presence
of gravity, we get a dynamics that is $SO(2)$-invariant but not $SO(3)$-invariant (i.e., it is invariant under a subgroup of the full symmetry group), contrary to what happens for the free rigid body \cite{HMR}. In the context of motion planning, symmetry breaking appears naturally in the form of a barrier function as we shall see in this work.

Symmetries in the condition of Pontryagin’s maximum principle (PMP) were studied in \cite{borum},  \cite{borumCDC}, \cite{eche}, \cite{manolo}, \cite{Tomoki-CDC}, \cite{Tomoki}. As an alternative to the PMP, a Lagrangian approach can be employed to exploit symmetries in the necessary conditions for optimal control problems. In the work \cite{BCGO}, Lagrangian systems with broken symmetries have been applied to optimal control problems on Lie groups for the single agent case. Also, in the single agent case, optimal control on semidirect products was also studied by using a Clebsch formulation \cite{GB}.  Both approaches provide coordinate-free equations describing necessary conditions for extrema on vector spaces.

Hence, in this work we do not solve differential equations on manifolds which usually requires embedding the manifold into an Euclidean space, and thus increasing the dimension of the system with additional constraints (differentiable manifolds are modeled locally on normed spaces and Lie groups are examples of differential manifolds). Alternatively, in this work we provide reduced equations on vector spaces which is useful in practice. In addition, in this work we build on the results of \cite{BCGO} by considering also obstacles in the configuration space that agents should avoid while they also avoid mutual collision.

%%%%%%%%%%%%%%%%%%%%%%%%%%%%%%%%%
The structure of the paper is as follows. Section II reviews the basic concepts of mechanics on manifolds, in particular, mechanics on Lie groups. Sections III describes left-invariant affine multi-agent control system and the problem addressed. Section IV presents the main results of the work, where we develop the reduction for optimality conditions in the optimal control problem. Section V illustrates the employment of the results with an example of three unicycles avoiding collisions and obstacles. Conclusions and future work are discussed in Section VI.

\section{Background}

\subsection{Mechanics on Manifolds}

Let $Q$ be a differentiable manifold, the configuration space of a mechanical system, with $\dim Q=n$. By $TQ$ we denote the tangent bundle of $Q$ and by $T^*Q$ its cotangent bundle. In local coordinates $Q$ is described by the positions $q=(q^1,\dots,q^n)$, a tangent vector $v_q=(q,\dot{q})\in TQ$ describes positions and the velocities $(q,\dot{q})=(q^1,\dots,q^n,\dot{q}^1,\dots,\dot{q}^n)$ and a cotangent vector $p_q=(q,p)=(q^1,\dots,q^n,p_1,\dots,p_n)\in T^*Q$ describes positions and the momenta.

\begin{definition}
Consider a differentiable function $f:Q\to M$ where $Q$ and $M$ are smooth manifolds. The \textit{tangent lift} (or tangent map) of $f$ at $q\in Q$ is the map $T_qf:T_qQ\to T_{f(q)}M$ and the \textit{cotangent lift} (or cotangent map) of $f$ is the map $T^*_qf:T^*_sM\to T^*_qQ$ defined by 
\begin{equation}\label{cotangent-lift}
    \langle T^*_qf(\alpha_s),v_q\rangle=\langle \alpha_s,T_qf(v_q)\rangle
\end{equation}

\noindent where $\alpha_s\in T^*_sM, \; v_q\in T_qQ, s=f(q)$ and $\langle\cdot,\cdot\rangle$ denotes the natural dual pairing between covectors and tangent vectors.

\end{definition}

The dynamics of a mechanical system is described by the equations of motion determined by a Lagrangian function $L:TQ\to\R$ given by $L(q,\dot{q})=K(q,\dot{q})-V(q)$, where $K:TQ\to\R$ denotes the kinetic energy and $V:Q\to\R$ the potential energy of the system. The equations of motion are given by the Euler-Lagrange equations 
$\displaystyle{\frac{d}{dt}\bigg(\frac{\partial L}{\partial\dot{q}^i}\bigg)-\frac{\partial L}{\partial q^i}=0, \; \; i=1,\dots,n}$, which determine a system of second-order differential equations. We assume the Lagrangian function is hyperregular and hence, through the Legendre transform $p=\frac{\partial L}{\partial\dot{q}}$, the Euler-Lagrange equations are equivalent to Hamilton's equations of motion $\displaystyle{\dot{q}=\frac{\partial H}{\partial p}, \; \; \dot{p}=-\frac{\partial H}{\partial q}}$, for the Hamiltonian $H(q,p)=\langle p,\dot{q}(q,p)\rangle - L(q,\dot{q}(q,p)),$ where $\dot{q}(q,p)$ is the second component of the inverse Legendre transform. Hamilton's equations can be defined equivalently in terms of the canonical Poisson bracket as follows.

\begin{definition} Let $P$ be a smooth manifold. A \textit{Poisson bracket} on $P$ is a bilinear, skew-symmetric operator $\{\cdot,\cdot\}:C^\infty(P)\times C^\infty(P)\to C^\infty(P)$ satisfying the Jacobi identity and the Leibniz rule. The pair $(P,\{\cdot,\cdot\})$ is called \textit{Poisson manifold}. For any Hamiltonian function $H:P\to\R$ the corresponding Hamiltonian vector field, $X_H$, which describes the equations of motion, is uniquely determined by $\dot{F}=\{F,H\}$ for all smooth $F:P\to\R$, along solutions to $X_H$ (see Def $4.22$ in \cite{HSS}).
\end{definition}

Consider the cotangent bundle of a manifold $Q$, for $F,G\in C^\infty(T^*Q)$ the pair $(T^*Q,\{\cdot,\cdot\})$ is a Poisson manifold with Poisson bracket $\{F,G\}=\sum_{i=1}^n\bigg(\frac{\partial F}{\partial q_i}\cdot\frac{\partial G}{\partial p_i}-\frac{\partial F}{\partial p_i}\cdot\frac{\partial G}{\partial q_i}\bigg)$. This  bracket is called the \textit{canonical Poisson bracket} (for more details see \cite{HSS}, \cite{MR}).

\subsection{Mechanics on Lie groups}

If the configuration space is a Lie group $G$, and the Lagrangian possesses the full $G$-symmetry, Euler-Lagrange equations can be reduced to a first-order system of differential equations.

\begin{definition}
Let $G$ be a Lie group and $Q$ a smooth manifold. A \textit{left-action} of $G$ on $Q$ is a smooth map $\Phi:G\times Q\to Q$ such that $\Phi(e,g)=g$ and $\Phi(h,\Phi(g,q))=\Phi(hg,q)$ for all $g,h\in G$ and $q\in Q$, where $e$ is the identity of the group $G$ and the map $\Phi_g:Q\to Q$ given by $\Phi_g(q)=\Phi(g,q)$ is a diffeomorphism for all $g\in G$.
\end{definition}

\begin{definition}
A function $f:Q\to\R$ is called \textit{invariant} under an action $\Phi_g$ if $f\circ\Phi_g=f$ for any $g\in G$.% 
\end{definition}

For a finite dimensional Lie group $G$, its Lie algebra $\g$ is defined as the tangent space to $G$ at the identity, $\g:=T_eG$. Let $L_g:G\to G$ be the left translation of the element $g\in G$ given by $L_g(h)=gh$ where $h\in G$. It is a left-action and a diffeomorphism on $G$. The tangent map of $L_g$ at $h\in G$ is $T_hL_g:T_hG\to T_{gh}G$ and the cotangent map is $T^*_hL_g:T^*_{gh}G\to T^*_hG$. By $\ad:\g\times\g\to\g,$ $(\xi,\eta)\mapsto\ad_\xi\eta$ we denote the adjoint operator which is defined by the Lie bracket on $T_eG$, i.e. $\ad_\xi\eta=[\xi,\eta]$, and by $\ad^*:\g\times\g^*\to\g^*$, $(\xi,\mu)\mapsto\ad^*_\xi\mu$ the coadjoint operator given by $\langle\ad^*_\xi,\mu,\eta\rangle=\langle\mu,\ad_\xi\eta\rangle$ where $\eta\in\g$ and $\langle\cdot,\cdot\rangle:\g^*\times\g\to\R$ is the  pairing between vectors and covectors (see \cite{V} and \cite{MR} for instance).

\begin{definition}
Denote by $\mathfrak{X}(G)$ the set of vector fields on $G$. A vector field $X:G\to TG,$ $h\mapsto X(h),$ is called \textit{left-invariant} if $T_hL_g(X(h))=X(L_g(h))$, for all $g,h\in G.$
\end{definition}

\begin{definition}
Consider a Lie group $G$, a vector space $W$ and the representation of $G$ on $W$ given by $\rho:G\times W\to W$, $(g,v)\mapsto\rho_g(v)$, which is a left action, and it is defined by the relation $\rho_{g_1}(\rho_{g_2}(v))=\rho_{g_1g_2}(v)$, $g_1,g_2\in G$. Its dual is given by $\rho^*:G\times W^*\to W^*,$ $(g,\alpha)\mapsto\rho^*_g(\alpha),$ satisfying $\langle\rho^*_{g^{-1}}(\alpha),v\rangle=\langle\alpha,\rho_{g^{-1}}(v)\rangle$.
\end{definition}

 The infinitesimal generator of the left action of $G$ on $W$ is $\rho':\g\times W\to W,$ $(\xi,v)\mapsto\rho'(\xi,v)=\frac{d}{dt}|_{t=0}\rho_{e^{t\xi}}(v).$ For every $v\in W$ consider the linear tranformation $\rho'_v:\g\to W,$ $\xi\mapsto\rho'_v(\xi)=\rho'(\xi,v)$ and its dual $\rho'^*_v: W^*\to\g^*,$ $\alpha\mapsto\rho'^*_v(\alpha).$ The last transformation defines the momentum map $\textbf{J}_W:W\times W^*\to\g^*,$ $(v,\alpha)\mapsto \textbf{J}_W(v,\alpha):=\rho'^*_v(\alpha)$ such that for every $\xi\in\g$, $\displaystyle{\langle \textbf{J}_W(v,\alpha),\xi\rangle=\langle\rho'^*_v(\alpha),\xi\rangle=\langle \alpha,\rho'_v(\xi)\rangle=\langle\alpha,\rho'(\xi,v)\rangle}$. 
 
 For $\xi\in\g$, consider the map $\rho'_\xi:W\to W$, $v\mapsto\rho'_\xi(v)=\rho'(\xi,v)$ and its dual $\rho'^*_\xi:W^*\to W^*$, $\alpha\mapsto\rho'^*_\xi(\alpha)$ such that $\langle\rho'^*_\xi(\alpha),v\rangle=\langle\alpha,\rho'_\xi(v)\rangle$, and hence it satisfies $\displaystyle{\langle\mathbf{J}_W(v,\alpha),\xi\rangle=\langle\alpha,\rho'(\xi,v)\rangle=\langle\alpha,\rho'_\xi(v)\rangle=\langle\rho'^*_\xi(\alpha),v\rangle}$. See \cite{HSS} and \cite{MR} for more details on the momentum map.

\section{Left-invariant multi-agent control systems}

Denote by $\mathcal{N}$ a set consisting of $s$ free agents, and by $\mathcal{E}\subset\mathcal{N}\times\mathcal{N}$ the set describing the interaction between agents. The neighbor relationships are described by an undirected graph $\mathcal{G}=(\mathcal{N},\mathcal{E})$, where $\mathcal{N}$ is the set of vertices and $\mathcal{E}$ the set of edges for $\mathcal{G}$. We further assume $\mathcal{G}$ is static and connected. We assume that each agent occupies a circular area of radius $\overline{r}$. For every agent $i\in \mathcal{N}$ the set $\mathcal{N}_i=\{j\in\mathcal{N}:(i,j)\in\mathcal{E}\}$ denotes the neighbors of that agent. The agent $i\in\mathcal{N}$ evolves on an $n$-dimensional Lie group $G$ and its configuration is denoted by $g_i\in G$. We denote by $\bf{G}$ and by $T_{\bf{e}}\bf{G}=:\bar{\g}$ the cartesian products of $s$ copies of $G$ and $\g$, respectively, where $\textbf{e}=\underbrace{(e,e,\dots,e)}_{s-copies}$ is the identity of $\bf{G}$ with $e$ the identity element of the of $G$.

For each agent $i\in\mathcal{N}$ there is an associated left-invariant control system described by the kinematic equations \begin{equation}\label{kin-each-agent}
    \dot{g}_i=T_{e}L_{g_i}(u_i), \;\; g_i(0)=g^0_i,
\end{equation} where $g_i(t)\in C^1([0,T],G), \; T\in\mathbb{R}$ fixed, $u_i\colon [0,T] \to \mathfrak{g}.$ is the control input and $g^i_0\in G$ is considered as the initial state condition. Note that while for each $i\in\mathcal{N}$, $\dim\g=n$ with $\g=\hbox{span}\{e_1,e_2,\ldots, e_n\}$, then the control inputs may be described by $u_i=[u_i^1,u_i^2,\dots,u_i^m]^T$, where $u_i(t)\in \mathbb{R}^{m}$, with $m\leq n$. Hence, the control input for each agent is given by $\displaystyle{u_i(t)=e_0+\sum_{k=1}^{m}u^k_i(t)e_k}$, where $e_0\in\g.$ Thus, the left-invariant control systems \eqref{kin-each-agent} for each agent $i\in\mathcal{N}$ can be written as \begin{equation}\label{kin-each-agent-basis}
    \dot{g}_i(t)=g_i(t)\bigg(e_0+\sum_{k=1}^{m}u^k_i(t)e_k\bigg).
\end{equation}

The problem under study consists on finding reduced necessary conditions for optimality in an optimal control problem for left-invariant multi-agent control system. These solution curves should minimize some cost function, prevent agent collision while also agents avoid static obstacles in the configuration space. 

\vspace{.2cm}

\textbf{Problem (collision and obstacle avoidance problem):} Find reduced optimality conditions on  $g(t)=(g_1(t),g_2(t),\dots,g_s(t))\in \textbf{G}$ and the controls $u(t)=(u_1(t),u_2(t),\dots,u_s(t))\in \bar{\g}$ such that  $(g(t),u(t))\in \textbf{G}\times\bar{\g}$ minimize the cost \begin{align}
	\min_{(g,u)}\sum_{i=1}^s\int_0^T\Big(C_i(g_i(t),u_i(t))&+U_i^{0}(g_i)\label{OCP}\\&+\frac{1}{2}\sum_{j\in \mathcal{N}_i}U_{ij}(g_i(t),g_j(t))\Big)dt \nonumber
\end{align}

\noindent subject to the kinematics $\dot{g}_i(t)=T_{e}L_{g_i(t)}(u_i(t))$, boundary values $g(0)=(g_1(0),\dots,g_s(0))=(g_1^0,\dots,g_s^0)$ and $g(T)=(g_1(T),\dots,g_s(T))=(g_1^T,\dots,g_s^T)$, under the following assumptions:

\begin{enumerate}[(i)]
\item{There is a left representation $\rho$ of $G$ on a vector space $W$, i.e., $\rho: G\to\mathrm{GL}(W)$ is a homomorphism.}
\item{$C_i: G\times\g\to\mathbb{R}$ are $G$-invariant functions for each  $i\in\mathcal{N}$ (under a suitable left action of $G$ on $G\times\g$, which will be defined shortly) and is also sufficiently regular.\label{a1}}
\item{$U_{ij}: G\times G\to\mathbb{R}$ (collision avoidance potential functions) are $G$-invariant functions under $\Phi$, defined  by \begin{align}
		\Phi: G \times (G \times G)&\longrightarrow G\times G,\label{eq_phi}\\
		(g,(g_1, g_2))&\longmapsto(L_{g}(g_1),L_{g}(g_2)),\nonumber
	\end{align} i.e., $U_{ij}\circ\Phi_{g} = U_{ij}$, for any $g \in G$, that is,  $U_{ij}(L_{g}(g_i),L_g(g_j)) = U_{ij}(g_i,g_j)$, for any $(g_i,g_j)\in \mathcal{E}$, $j\in\mathcal{N}_i$ and they are also sufficiently regular.}

\item{$U_{i}^{0}: G\to\mathbb{R}$ (obstacle avoidance potential functions) are not $G$-invariant functions and they are also sufficiently regular, for $i\in\mathcal{N}$.}

\item{The obstacle avoidance functions $U_{i}^{0}$ depend on a parameter $\alpha_{0}\in W^{*}$ for each agent $i\in\mathcal{N}$. Hence, we can define the extended potential function as $U_{i,0}^{\textnormal{ext}}: G\times W^{*}\to\mathbb{R}$, with $\Vezero(\cdot,\alpha_{0}) = U_{i}^{0}$.}

\item{The extended obstacle avoidance functions are $G$-invariant under 
	$\tilde{\Phi}$, defined by \begin{align}
		\tilde{\Phi}: G \times (G \times W^{*})&\longrightarrow G\times W^{*},\label{eq_phi_tilde}\\
		(g,(h,\alpha))&\longmapsto(L_{g}(h),\rho_{g^{-1}}^{*}(\alpha)),\nonumber
	\end{align}  where $\rho_{g^{-1}}^{*}\in\mathrm{GL}(W^{*})$ is the adjoint of $\rho_{g^{-1}}\in\mathrm{GL}(W)$, i.e., $U_{i,0}^{\textnormal{ext}}\circ\tilde{\Phi}_{g} = U_{i,0}^{\textnormal{ext}}$, for any $g \in G$, or $U_{i,0}^{\textnormal{ext}}(L_{g}(h),\rho_{g^{-1}}^{*}(\alpha)) = U_{i,0}^{\textnormal{ext}}(g,\alpha)$ where $\alpha\in W^{*}$}.

\item{The obstacle avoidance potential functions are invariant under the left action of the isotropy group
	\begin{equation}
		G_{\alpha_{0}}=\{g\in G\mid\rho_{g}^{*}(\alpha_{0}) = \alpha_{0}\}.
	\end{equation}\label{a6}}

\end{enumerate}

\vspace{-0.2cm}

 Note that $G\times\g$ is a trivial vector bundle over $G$ and define the left action of $G$ on $G\times\g$ as follows
\begin{align}
\Psi: G\times (G\times\g)&\longrightarrow (G\times\g),\nonumber\\
(g,(h,u))&\longmapsto(L_{g}(h),u).\label{eq_psi}
\end{align}

We further assume that $C_i: G\times\g\to\mathbb{R}$ is $G$-invariant under \eqref{eq_psi}, i.e., $C_i\circ\Psi_{g} = C_i$, for any $g\in G$. 

\section{Reduction of necessary conditions}

We can solve the proposed problem as a constrained problem using the method of Lagrange multipliers (see, e.g, \cite{B}, \cite{BCGO},  \cite{KM}) or by using the Pontryagin Maxiumum principle (PMP) \cite{K}, \cite{JK2}, \cite{Tomoki}. Our approach is based on the PMP. 

Define the \textit{augmented Hamiltonian} function $H_{a}: \textbf{G}\times\bar{\mathfrak{g}}\times\bar{\mathfrak{g}}^{*}\times (W^{*})^s\to\mathbb{R}$ as follows
\begin{equation*}
    \begin{split}
         H_{a}(g,u,\mu,\alpha)=&\sum_{i=1}^s\bigg[\langle\mu_i,u_i-e_{0}\rangle-C_i(g_i,u_i-e_{0})\\
        &-\Vezero(g_i,\alpha)-\frac{1}{2}\sum_{j\in\mathcal{N}_i}U_{ij}(g_i,g_j)\bigg],
    \end{split}
\end{equation*}
where $\mu(\cdot)\in C^{1}([0,T],\mathfrak{g}^{*})$. By applying all actions of the previous section component-wise each component of the augmented Hamiltonian is $G$-invariant thus we define the \textit{reduced augmented Hamiltonian} $h_{a}:G^{s-1}\times\bar{\g}\times\bar{\mathfrak{g}}^{*}\times (W^{*})^s\to\mathbb{R}$, which is given by
\begin{equation*}
    \begin{split}
         h_{a}(g,u,\mu,\alpha)=&\sum_{i=1}^s\bigg[\langle\mu_i,u_i-e_{0}\rangle-C_i(u_i-e_{0})\\
        &-\Vezero(\alpha_i)-\frac{1}{2}\sum_{j\in\mathcal{N}_i}U_{ij}(g_i^{-1}g_j)\bigg],
    \end{split}
\end{equation*}
where $\alpha_i = \rho^{*}_{g_i}(\alpha)$ and with a slight abuse of notation, we write $C_i(e,u_i-e_{0}) = C_i(u_i-e_{0})$, $\Vezero(\alpha_i) = \Vezero(e,\alpha_i)$ and $U_{ij}(e,g_i^{-1}g_j)=U_{ij}(g_i^{-1}g_j)$. 

Next, we can employ the Pontryagin Maximum Principle for the reduced augmented Hamiltonian $h_a$. Therefore, it follows that necessary conditions for normal extrema for the OCP can be obtained by the integral curves of the Hamiltonian vector field for the \textit{reduced optimal Hamiltonian} $h:G^{s-1}\times\bar{\mathfrak{g}}^{*}\times (W^{*})^s\to\mathbb{R}$ given by $$h(g,\mu,\alpha)=\max_{u}h_{a}(g,u,\mu,\alpha)=h_a(g,u^{*},\mu,\alpha),$$where $u^{*}$ denotes the optimal control. The reduced optimal Hamiltonian $h:G^{s-1}\times\bar{\mathfrak{g}}^{*}\times (W^{*})^s\to\mathbb{R}$ can now be obtained and is given by 
\begin{equation*}
    \begin{split}
         h(g,\mu,\alpha)=&\sum_{i=1}^s\bigg[\langle\mu_i,u^*_i-e_0\rangle-C_i(u^*_i-e_0)-\Vezero(\alpha_i)\\
        &-\frac{1}{2}\sum_{j\in\mathcal{N}_i}U_{ij}(g_i^{-1}g_j)\bigg],
    \end{split}
\end{equation*}
where $\alpha_i = \rho^{*}_{g_i}(\alpha)$.

Hamilton's equations of motion are determined from the equation $\dot{F}=\{F,h\}$ where $\{\cdot,\cdot\}$ is the Lie-Poisson bracket on $G^{s-1}\times\bar{\g}^*\times (W^*)^s$ and $F\in\mathcal{
C}^{\infty}(G^{s-1}\times\bar{\g}^*\times (W^*)^s)$. The left hand side of Hamilton's equations is $\displaystyle{\dot{F}=\frac{dF}{dt}=DF\cdot(\dot{g},\dot{\mu},\dot{\alpha})=\Big\langle(\dot{g},\dot{\mu},\dot{\alpha}), \frac{\partial F}{\partial(g,\mu,\alpha)}\Big\rangle}$. By Prop. 4.4.1 of \cite{AM}, the Lie-Poisson bracket on $G^{s-1}\times\bar{\g}^*$ is given by \begin{equation*}
\begin{split}
	\{F,h\}(g,\mu)=&\sum_{i=1}^{s}\Bigg[\Big\langle T_e^*L_{g_i}\Big(\frac{\partial F}{\partial g_i}\Big), \frac{\partial h}{\partial\mu_i}\Big\rangle \\ 
	&+ \Big\langle T_e^*L_{g_i}\Big(\frac{\partial h}{\partial g_i}\Big),\frac{\partial F}{\partial\mu_i}\Big\rangle
	- \Big\langle\mu_i,\Big[\frac{\partial F}{\partial\mu_i},\frac{\partial h}{\partial\mu_i}\Big]\Big\rangle\Bigg],
\end{split}	
\end{equation*}

\noindent where $\frac{\partial F}{\partial g_i}$ and $\frac{\partial h}{\partial g_i}$ are seen as covectors and $[\cdot,\cdot]$ is the Lie bracket on the Lie algebra $\g$. We apply this Lie-Poisson bracket in our setting hence, so the Lie-Poisson bracket on $G^{s-1}\times\bar{\g}^*\times (W^*)^s$ is given by \begin{equation}\label{Lie-Poisson brac}
\begin{split}
\{F,h\}(g,\mu,\alpha)=&\sum_{i=1}^{s}\Bigg[\Big\langle T_e^*L_{g_i}\Big(\frac{\partial F}{\partial g_i}\Big), \frac{\partial h}{\partial\mu_i}\Big\rangle \\ 
+& \Big\langle T_e^*L_{g_i}\Big(\frac{\partial h}{\partial g_i}\Big),\frac{\partial F}{\partial\mu_i}\Big\rangle- \Big\langle\mu_i,\Big[\frac{\partial F}{\partial\mu_i},\frac{\partial h}{\partial\mu_i}\Big]\Big\rangle \\
-& \Big\langle\alpha_i,\frac{\partial F}{\partial\mu_i}\frac{\partial h}{\partial\alpha_i}-\frac{\partial h}{\partial\mu_i}\frac{\partial F}{\partial\alpha_i}\Big\rangle\Bigg],
\end{split}	
\end{equation} (for more details one can see Section $9.6$ of \cite{HSS}). 

From the first term of the previous equation, through the property of the cotangent lift \eqref{cotangent-lift}, we get:

\[\sum_{i=1}^s\Big\langle T_e^*L_{g_i}\Big(\frac{\partial F}{\partial g_i}\Big), \frac{\partial h}{\partial\mu_i}\Big\rangle=\sum_{i=1}^s\Big\langle T_eL_{g_i}\Big(\frac{\partial h}{\partial\mu_i}\Big),\frac{\partial F}{\partial g_i}\Big\rangle.\]

The Lie bracket in equation (\ref{Lie-Poisson brac}) is the adjoint map, thus each component of the sum of the third term yields

\begin{equation*}
\begin{split}
\Big\langle\mu_i,\Big[\frac{\partial F}{\partial\mu_i},\frac{\partial h}{\partial\mu_i}\Big]\Big\rangle &= \Big\langle\mu_i,\ad_\frac{\partial F}{\partial\mu_i}\frac{\partial h}{\partial\mu_i}\Big\rangle = -\Big\langle\mu_i,\ad_\frac{\partial h}{\partial\mu_i}\frac{\partial F}{\partial\mu_i}\Big\rangle \\
&= -\Big\langle\ad^*_\frac{\partial h}{\partial\mu_i}\mu_i,\frac{\partial F}{\partial\mu_i}\Big\rangle.
\end{split}
\end{equation*}

\noindent The last term can be written as

\begin{multline*}
	\sum_{i=1}^s\Bigg[\Big\langle\alpha_i,\frac{\partial F}{\partial\mu_i}\frac{\partial h}{\partial\alpha_i}\Big\rangle-\Big\langle\alpha_i,\frac{\partial h}{\partial\mu_i}\frac{\partial F}{\partial\alpha_i}\Big\rangle\Bigg] = \\ =\sum_{i=1}^s\Bigg[\Big\langle\mathbf{J}_W\Big(\frac{\partial h}{\partial\alpha_i},\alpha_i\Big), \frac{\partial F}{\partial\mu_i}\Big\rangle + \Big\langle\frac{\partial h}{\partial\mu_i}\alpha_i,\frac{\partial F}{\partial\alpha_i}\Big\rangle\Bigg],
\end{multline*} where $\mathbf{J}_W$ is the momentum map and the first entry of the last term is the dual of the infinitesimal generator of the left action of $G$ on $W$. Thus, for the reduced Hamiltonian function Lie-Poisson equations of motion, for each agent, are given by:

\begin{eqnarray}\label{Lie-Poiss eqs}
&&\dot{g_i}=T_{e}L_{g_i}(u_i^*-e_0), \nonumber \\
&&\dot{\mu_i}=\ad^*_{u_i}\mu_i +\sum_{j\in\mathcal{N}_i}T^*_eL_{g_i}\Big(\frac{\partial U_{ij}}{\partial g_i}\Big) -\mathbf{J}_W\Big(\frac{\partial \Vezero}{\partial\alpha_i},\alpha_i\Big), \nonumber \\
&&\dot{\alpha_i}=-\rho'^*_{u_i}(\alpha_i).
\end{eqnarray}

\section{Application to optimal control of unicycles}

A unicycle is a homogeneous disk rolling vertically on a horizontal plane. Consider $s$ such unicycles. The configuration space whose elements determine the motion of each unicycle is $SE(2)\cong SO(2)\times\R$. An element $g_i\in SE(2)$ is given by $g_i=\begin{pmatrix}
\cos\theta_i&-\sin\theta_i&x_i\\
\sin\theta_i&\cos\theta_i&y_i\\
0          &          0& 1
\end{pmatrix}$, where $(x_i,y_i)\in\R^2$ represents the point of contact of the disk with the horizontal plane and $\theta_i$ represents the angular orientation of the unicycle. The control inputs, for each agent, are given by $u_i=(u_i^1,u_i^2)$ where $u_i^1$ denotes the steering speed and $u_i^2$ denotes the rolling speed of the unicycle. For more details see \cite{B}. 

The kinematic equations for the multi-agent system are:
\begin{equation}\label{kinem-eqs}
  \dot{x_i}=u_i^2\cos\theta_i, \; \; \dot{y_i}=u_i^2\sin\theta_i,  \; \; \dot{\theta_i}=u_i^1,\, i=1,\ldots,s.  
\end{equation}

The Lie algebra $\se$ of $SE(2)$ is described by $\se=\Big\{\begin{pmatrix}
A&b\\
0&0
\end{pmatrix}: A\in\mathfrak{so}(2)\text{ and }b\in\R^2\Big\}$ where $A=-aJ$ with $a\in\R$ and $J=\begin{pmatrix}
0&1\\
-1&0
\end{pmatrix}$. The Lie algebra $\se$ is identified with $\R^2$ through the isomorphism $\begin{pmatrix}
-aJ&b\\
0&0
\end{pmatrix}\mapsto (a,b).$
The elements of the basis of the Lie algebra $\se$ are 

\begin{equation*}
e_1=\begin{pmatrix}
0&-1&0\\
1&0&0\\
0&0&0
\end{pmatrix}, 
e_2=\begin{pmatrix}
0&0&1\\
0&0&0\\
0&0&0
\end{pmatrix},
e_3=\begin{pmatrix}
0&0&0\\
0&0&1\\
0&0&0
\end{pmatrix},
\end{equation*}

\noindent which satisfy $[e_1,e_2]=e_3,\; [e_2,e_3]=0, \; [e_3,e_1]=e_2$. Thus, the kinematic equations (\ref{kinem-eqs}) take the form $\dot{g_i}=g_iu_i=g_i(u_i^1e_1+u_i^2e_2)$ and give rise to a left-invariant control system on $SE(2)^s\times \mathfrak{se}(2)^s$.

The inner product on $\se$ is given by $\langle\xi_1,\xi_2\rangle=\tr(\xi_1^T\xi_2)$ for $\xi_1,\xi_2\in\se$ and hence, the norm is given by $\|\xi\|=\sqrt{\tr(\xi^T\xi)},$ for any $\xi\in\se$.

The dual Lie algebra $\se^*$ of $SE(2)$ is defined through the dual pairing, $\langle\alpha,\xi\rangle=\tr(\alpha\xi)$, where $\alpha\in\se^*$ and $\xi\in\se$ hence, the elements of the basis of $\se^*$ are \begin{equation*}
e^1=\begin{pmatrix}
0&\frac{1}{2}&0\\
-\frac{1}{2}&0&0\\
0&0&0
\end{pmatrix}, 
e^2=\begin{pmatrix}
0&0&0\\
0&0&0\\
1&0&0
\end{pmatrix},
e^3=\begin{pmatrix}
0&0&0\\
0&0&0\\
0&1&0
\end{pmatrix}.
\end{equation*}

Consider the cost function $C_i(g_i,u_i)=\frac{1}{2}\langle u_i,u_i\rangle$ and the artificial potential function $U_{ij}:SE(2)\times SE(2)\to\R$ given by $\displaystyle{U_{ij}(g_i,g_j)=\frac{\sigma_{ij}}{2((x_i-x_j)^2+(y_i-y_j)^2-4\overline{r}^2)}}$, where $\sigma_{ij}\geq 0$ and $\overline{r}$ is the radius of the disk each agent occupies as defined at the end of Section III. Consider a spherical obstacle with unit radius and without loss of generality let it be centered at the origin. Hence, consider the obstacle avoidance potential function $U_i^0:SE(2)\to\R$, $\displaystyle{U_i^0(g_i)=\frac{\sigma_{i0}}{2(x^2+y^2-(\overline{r}+1)^2)}}$, where $\sigma_{i0}>0$. 

Note that the obstacle avoidance potential functions are not $SE(2)$-invariant but $SO(2)$-invariant, so they break the symmetry. Using the norm of $\se$ and for $\overline{r}=1$, $U_{ij}$ and $U_i^0$ are equivalently given by
$\displaystyle{U_{ij}(g_i,g_j)=\frac{\sigma_{ij}}{2(\|\Ad_{g_i^{-1}g_j}e_1\|^2-6)}}$ and 
$\displaystyle{U_i^0(g_i)=\frac{\sigma_{i0}}{2(\|\Ad_{g^{-1}}e_1\|^2-6)}}$.

Let $W=\g^*$, so we define the extended potential functions $\Vezero:SE(2)\times \g\to\R$ by
$\Vezero(g_i,\alpha)=\frac{\sigma_{i0}}{2(\|\Ad_{g^{-1}}\alpha\|^2-6)}$, which are $SE(2)$-invariant under the action of $\tilde{\Phi}$, i.e. $\Vezero\circ\tilde{\Phi}=\Vezero,$ for any $g\in SE(2)$. Since, $W=\g^*$ we have $\displaystyle{\mathbf{J}_W\Big(\frac{\partial\Vezero}{\partial\alpha_i},\alpha_i\Big)=\ad^*_{\alpha_i}\Big(\frac{\partial\Vezero}{\partial\alpha_i}\Big)}$, $\rho'^*_{u_i}(\alpha_i)=\ad_{u_i}\alpha_i,$ and $\rho^*_{g_i}(\alpha)=\Ad_{g^{-1}}\alpha$ (for more details see Corollary 7 of \cite{BCGO}), recall from section III that $\rho^*$ is the adjoint of $\rho$. Thus, from PMP we get $u_i^1=\frac{1}{2}\mu_i^1$, $u_i^2=\mu_i^2$ and $u_3^i=\mu_i^3$ and the reduced augmented Hamiltonian $h:G^{s-1}\times(\g^*)^s\times W^* \to \R$ is given by \begin{equation*}
	\begin{split}
		h(g,\mu,\alpha)=&-\frac{3}{4}\sum_{i=1}^{s}(\mu_i^1)^2+\frac{1}{2}\sum_{i=1}^{s}(\mu_i^2)^2\\ &+\Vezero(e,\alpha_i)+\frac{1}{2}\sum_{j\in \mathcal{N}_i}U_{ij}(e,g_i^{-1}(t)g_j(t)).
	\end{split}
\end{equation*}  Denote by $\Gamma_{ij}:=T_{\overline{e}_i}^{*}L_{g_i}\left(\frac{\partial U_{ij}}{\partial g_i}\right)$ and note that \begin{small}\begin{align*}
    \ad^*_{u_i}\mu_i &=\begin{pmatrix}
0&-\frac{\mu_i^2\mu_i^3}{2}&0\\
\frac{\mu_i^2\mu_i^3}{2}&0&0\\
\frac{\mu_i^1\mu_i^3}{2}&-\frac{\mu_i^1\mu_i^2}{2}&0
\end{pmatrix},   \\
 & \\
\ad^*_{\alpha_i}\Big(\frac{\partial\Vezero}{\partial\alpha_i}\Big) &= \begin{pmatrix}
		0&0&0\\
		0&0&0\\
		\Gamma_{i,0}^{31}&\Gamma_{i,0}^{32}&0
		\end{pmatrix}  \\
 &=\frac{\sigma_{i0}\alpha_i^1}{(\|\alpha_i\|^2-6)^2}\begin{pmatrix}
		0&0&0\\
		0&0&0\\
		-\alpha_i^3&\alpha_i^2&0
		\end{pmatrix}, \\
 & \\
T^*_{\overline{e}_i}L_{g_i}\Big(\frac{\partial U_{ij}}{\partial g_i}\Big) &= \begin{pmatrix}
0&0&0\\
0&0&0\\
\Gamma_{ij}^{31}&\Gamma_{ij}^{32}&0
\end{pmatrix} \\
 &= \frac{-\sigma_{ij}}{((x_{ij})^2+(y_{ij})^2-4)^2}\begin{pmatrix}
0&0&0\\
0&0&0\\
x_{ij}&y_{ij}&0
\end{pmatrix},\\
 & \\
\ad_{u_i}\alpha_i &= \begin{pmatrix}
0&0&-u_i^1\alpha_i^3\\
0&0&u_i^1\alpha_i^2-u_i^2\alpha_i^1\\
0&0&0
\end{pmatrix},
\end{align*}where $x_{ij}=x_i-x_j, \; y_{ij}=y_i-y_j$ and

\[\Ad_{g_i^{-1}}\alpha_0=\begin{pmatrix}
0&-1&x_i\sin\theta_i-y_i\cos\theta_i\\
1&0&x_i\cos\theta_i+y_i\sin\theta_i\\
0&0&0
\end{pmatrix}.\] \end{small}

Thus, the Lie-Poisson equations (\ref{Lie-Poiss eqs}) are 

$\begin{array}{ll}
     & \\
    \dot{\mu}_i^1=-\frac{1}{2}\mu_i^2\mu_i^3, &  \\
     & \\
    \dot{\mu}_i^2=\frac{1}{2}\mu_i^1\mu_i^3-\Gamma_{i,0}^{31}+\sum_{j\in\mathcal{N}_i}\Gamma_{ij}^{31}, &  \\
     & \\
    \dot{\mu}_i^3=-\frac{1}{2}\mu_i^1\mu_i^2-\Gamma_{i,0}^{32}+\sum_{j\in\mathcal{N}_i}\Gamma_{ij}^{32}, & \\
    & 
\end{array}$

\noindent together with 
\vspace{0.3cm}

$\begin{array}{ll}
\dot{\alpha}_i^1=0, & \alpha_i^1=1, \\
\dot{\alpha}_i^2=\frac{1}{2}\mu_i^1\alpha_i^3,  &  \alpha_i^2=x_i\sin\theta_i-y_i\cos\theta_i,  \\
\dot{\alpha}_i^3=-\frac{1}{2}\mu_i^1\alpha_i^2+\mu_i^2\alpha_i^1, & \alpha_i^3=x_i\cos\theta_i+y_i\sin\theta_i. 
\end{array}$

\vspace{0.5cm}

\begin{remark}
	Note that, for the unicycle model, alternatively it can be used different potential functions with both properties (collision and obstacle avoidance). For instance one can consider the kind of potential functions $\displaystyle{U(g_i,g_j)=\frac{\tilde{\sigma}_i}{2(x_i^2+y_i^2-(\overline{r}+1)^2)\Upsilon_{ij}}}$, where $\Upsilon_{ij}=\prod_{j\in\mathcal{N}_i}((x_i-x_j)^2+(y_i-y_j)^2-d^2_{ij})$ and $d_{ij}$ stands for the desired distance between two agents which were used for instance in \cite{O-SM}.\hfill$\diamond$\end{remark}

Next we present some simulation results for the previous optimal control problem for three agents on $SE(2)$ avoiding mutual collision as well as a static obstacle in the configurations space. For simulation results, we transform the two-point boundary value problem to an initial value problem by considering the shooting method. The controlled dynamics of the three agents are simulated using Euler's method with time step $h=0.0001$ and $N=50000$. Initial conditions are $g_1(0)=\begin{pmatrix}
	1 & 0 & -1 \\
	0 & 1 & -5 \\
	0 & 0 & 1 
\end{pmatrix}$, 
$g_2(0)=\begin{pmatrix}
	\sqrt{2}/2 & \sqrt{2}/2 & -1 \\
	\sqrt{2}/2 & \sqrt{2}/2 & -7 \\
	0          &  0         & 1  
\end{pmatrix}$, 
$g_3(0)=\begin{pmatrix}
	1 & 0 & -1 \\
	0 & 1 &  0 \\
	0 & 0 & 1 
\end{pmatrix}$, with initial angular velocities $u^1_1(0)= u^1_2(0)=1$, $u^1_3(0)=1.5$ and linear velocities $u^2_1(0)=u^2_2(0)=0.6$, $u^2_3(0)=2$. 

\begin{figure}[h]
	\includegraphics[scale=0.3]{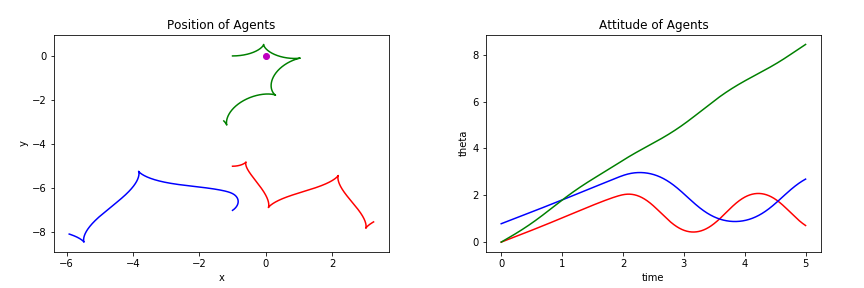}
	\caption{The left plot illustrates the trajectories of the three unicycles. The right plot shows evolution of their attitude through time. }
	\label{pos-att}
\end{figure}

\begin{figure}[h]
	\includegraphics[scale=0.3]{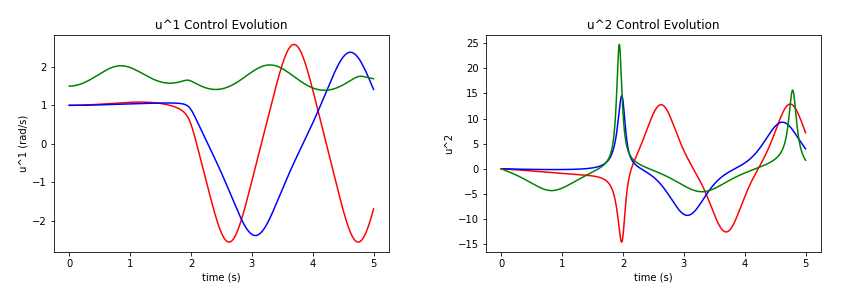}
	\caption{The two plots show the change of angular and linear velocities with respect to time for the three agents.}
	\label{u12}
\end{figure}

Fig.\ref{pos-att} shows the trajectories of the agents on the configuration space. Initial values were chosen in order that one of the agents avoids the obstacle and the other two prevent mutual collision. Fig.\ref{u12} shows the evolution of controls inputs through the time. Agents $1$, $2$ and $3$ are depicted by the color red, blue and green, respectively. Note that in Fig. \ref{u12} we can also observe that the agents (blue and red) avoiding mutual collision have similar angular velocities, while the agent avoiding the obstacle and depicted by green color, has less sharp fluctuations. From the second graph we can also see that almost at the same moment an increase of the linear velocities occurs so that all three agents prevent collision.

\section{Conclusions}

We studied the reduction by symmetries for necesary optimality conditions for extrema in an optimal control problem of left-invariant affine multi-agent control system, by exploiting the physical symmetries of the agents and static obstacles in the configuration space. Reduced optimality conditions were obtained by using techniques from Hamiltonian mechanics on Lie groups together with the PMP. We applied the results to a collision avoidance problem for multiple unicycles in the presence of an obstacle.

\end{document}